\documentclass[12pt,a4]{article}

\topmargin= - 10 true mm
\usepackage{amsmath,graphicx,amsthm,amssymb}
\usepackage[font=footnotesize,labelfont=bf]{caption}
\usepackage[numbers,square,sort&compress]{natbib}

\newcommand{\exn}{{\bf E}}
\newcommand{\pr}{{\bf P}}
\newcommand{\N}{\mathbb{N}}
\newcommand{\R}{\mathbb{R}}

\newcommand{\deq}{\stackrel{d}{=}}
\newcommand{\ep}{\varepsilon}

\newcommand{\Va}{\mbox{\rm Var}\, }
\newcommand{\Co}{\mbox{\rm Cov}\, }
\newcommand{\ind}{{\bf 1}}

\newcommand{\topr}{\stackrel{P}{\longrightarrow}}
\newcommand{\todi}{\stackrel{d}{\longrightarrow}}

\newcommand{\gest}{\stackrel{st}{\ge}}
\newcommand{\lest}{\stackrel{st}{\le}}

\newtheorem*{theo}{Theorem}

\newtheorem*{lemo}{Lemma}

\title{On the maximum of discretely sampled fractional Brownian motion with small Hurst parameter}
\author{Konstantin Borovkov$^1$ and Mikhail Zhitlukhin$^2$}

\date{}

\begin{document}
\maketitle

\footnotetext[1]{School of Mathematics and Statistics, The University of Melbourne, Parkville 3010, Australia; e-mail: borovkov@unimelb.edu.au.}

\footnotetext[2]{Steklov Mathematical Institute of Russian Academy of Sciences, Gubkina str.~8,
119991, Moscow, Russia; email: mikhailzh@mi.ras.ru.}

\begin{abstract}
We show that the distribution of the maximum of the fractional Brownian motion
$B^H$ with Hurst parameter $H\to 0$ over an $n$-point  set $\tau \subset [0,1]$ can be approximated  by the normal law  with mean $\sqrt{\ln n}$ and variance $1/2$ provided that $n\to \infty$ slowly enough and the points in~$\tau$ are not too close to each other.


\end{abstract}

\section{Introduction}

Let $\{B^H_t\}_{t\ge 0}$ be the fractional Brownian motion (fBM) with Hurst
index $H\in (0,1]$. Recall that the fBM is a zero-mean continuous Gaussian
process with the covariance function
\[
\exn B^H_s B^H_t=\frac12 \bigl(s^{2H} + t^{2H} - |t-s|^{2H}\bigr), \quad s,t\ge 0.
\]
Alternatively, $B^H$ can be defined as a continuous Gaussian process
with stationary increments such that $B_t^H$ has zero mean and variance~$t^{2H}$. In particular, $W:=B^{1/2}$ is the
standard Brownian motion (BM)  that has independent increments. The increments of $B^H$ are positively correlated if $H>1/2$ and negatively correlated if
$H<1/2$.

The fBM has found use in many models in applied fields (see, e.g., the survey in the preface to the monograph~\cite{M07}).
In particular, the processes~$B^H$ with small~$H$ (the 
case we are focussing on in this paper) have  recently  been used to model stock  price volatility~\cite{GJR,BFG}. It is interesting and important for a number of applications  to know the distribution (or a suitable approximation thereof) of the maximum
\[
\overline{B^H_T}:=\max_{0\le t\le T} B^H_t 
\]
of the fBM on a fixed time interval $[0,T]$, $T>0.$ Unfortunately, besides the case of the standard  BM ($H=1/2$) and the degenerate 
case $H=1$ (where $B_t^1 = \zeta t,$   $t\ge 0,$ for  a standard normal
random variable $\zeta$), there is no known closed form expression for the distribution of $\overline{B_T^H}.$  As in practice one usually deals with discretely sampled data, what would be of real practical interest is actually the behavior of the distribution of the maximum of the fBM  sampled on a discrete time grid on~$[0,T]$. 

In this paper, we consider the case when $H$ vanishes and deal with
the maxima of the fBM $B^H$  sampled on a
(generally speaking, non-uniform) discrete time grid.   
Recall that in that case the finite-dimensional distributions of $B^H$ converge to those of a ``translated'' continuum of independent normal random variables (see, e.g.,~\cite{BoMiNoZh}):
\begin{equation}
\{B^H\}_{t\ge 0 }\stackrel{\text{f.d.d.}}{\longrightarrow}  
  \{\xi\}_{t\ge 0} \quad \mbox{as} \quad   H\to 0, 
\label{Hto0}
\end{equation}
where $\xi_t := (\zeta_t-\zeta_0)/\sqrt{2},$   $\{\zeta_t\}_{t\ge 0}$ is 
a family of independent standard normal random variables. It is clear from~\eqref{Hto0} that $\overline{B^H_T} \topr \infty$ as
$H\to 0$. However, for  any fixed finite subset  
\begin{equation}
\tau=\{t_i\}_{i=1}^n \subset [0,T],\quad \mbox{where} \  t_1<t_2<\cdots<t_n ,
\label{tau}
\end{equation}
if one considers the random vector
\[
B^{H,\tau}:= (B^H_{t_1},B^H_{t_2},\ldots, B^H_{t_n} )\in\R^n,
\] 
and let 
$\overline{x} := \max_{1\le i\le n} x_i$ for a vector $x\in\R^n,$ 
relation~\eqref{Hto0} implies the convergence in distribution
\begin{equation}
\overline{B^{H,\tau}} \todi  (\overline{\zeta^n} -\zeta_0)/\sqrt{2} \quad \mbox{as} \quad H\to 0,
\label{Hnto0}
\end{equation}
where $\zeta^n:= (\zeta_1, \ldots, \zeta_n)$. One can
easily see that the distribution function of the random variable on the RHS
of~\eqref{Hnto0} is given by the convolution $(\Phi^n * \Phi) (\sqrt{2}x)$,
where $\Phi$ is the standard normal distribution function.

Now what can be said about the behavior of $\overline{B^{H,\tau}}$ when
simultaneously $H\to 0$ and the number~$n$ of points in the partition~$\tau$
tends to infinity? One can conjecture that, if $n\to\infty$ slowly  enough (so that the
dependence between the components of the vector  $B^{H,\tau}$ decays sufficiently quickly), then the distribution of
$\overline{B^{H,\tau}}$ would still be close to that of the RHS
of~\eqref{Hnto0}. The behavior of the distribution of $\overline{\zeta^n}$
as $n\to\infty$ has been  known since the work of Fisher and Tippett~\cite{FiTi}
who demonstrated that, taking $a_n:= \sqrt{2\ln n }$ and
$b_n:=\sqrt{2\ln n } - (\ln\ln n +\ln (4\pi))/(2\sqrt{2\ln n }),$ one has
\begin{equation}
 a_n(\overline{\zeta^n}   -b_n)\todi G \quad\mbox{as}\quad  n\to\infty,
\label{toGu}
\end{equation}
where the limiting random variable $G$ follows the Gumbel distribution
$\Lambda(x)=e^{-e^{-x}},$ $x\in\R$. In fact, the uniform distance between
the distribution functions of the LHS of~\eqref{toGu} and $\Lambda$ was
shown to be of the order of $1/\ln n$~\cite{Ha}. Choosing slightly different
sequences 
\begin{equation}
b_n:= \Phi^{-1}(1-1/n),\quad a_n:=b_n + 1/b_n,
\label{ab}
\end{equation}
one can show that
that distance admits an asymptotic upper bound of the form $1/(3 \ln n)$
(see~\cite{GaJoUt}).

So one can expect a first order approximation of the form
$\sqrt{\ln n}+\zeta_0/\sqrt{2}$ to hold true for the
maximum~$\overline{B^{H,\tau}}$ as $n\to\infty$, provided that $H\to 0$ fast
enough for the given decay rate of the distance between the points~$t_i$. Our
main result below  confirms that conjecture and
specifies conditions under which it holds. Without loss of generality, we 
consider  the case $T=1$ only, since the case of arbitrary~$T$ can be easily
reduced to the former using the self-similarity property of the fBM.

\section{The main result}

Denote by $\lest $ the stochastic order relation for random variables: we
write $\xi\lest \eta$ iff $\pr (\xi \le x) \ge \pr (\eta \le x)$, $x\in \R$,
and   $\xi\gest \eta$ iff $\eta\lest \xi$. By 
\[
\delta (\tau) := \min_{1\le i \le n } (t_{ i} - t_{i-1}) , \quad \mbox{where $t_{0}:=0$},
\]
we  denote the minimal distance between the points of the finite subset  $\tau$ (cf.~\eqref{tau}). 
As usual,  $o_P(1)$ denotes 
a sequence of random variables converging to zero in probability.

\begin{theo}
Let  $H_k\in(0,1]$ be such that $H_k\to 0$ as $k\to\infty$, and
$\tau_k=\{t_{k,i}\}_{i=1}^{n_k} $  be  a sequence of subsets of $(0,1]$, $ t_{k,1}<\cdots<t_{k,n_k}, $ such that $n_k\to\infty$, $\delta_k:=\delta (\tau_k).$

\smallskip

{\rm (i)} If $H_k(\ln n_k)^{1/2}\to 0$ and $H_k\ln (n_k\delta_k) \to 0$ as $k\to\infty$ then
\begin{equation}
\overline{B^{H_k,\tau_k}}\lest \sqrt{\ln n_k}+\zeta_0/\sqrt{2}  +o_P(1).
\label{<=}
\end{equation}

\smallskip

{\rm (ii)} If $H_k(\ln n_k)^{2}\to0$ and $H_k \ln\delta_k\to 0$ as
$k\to\infty$, then
\[
\overline{B^{H_k,\tau_k}}\gest \sqrt{\ln n_k} +\zeta_0/\sqrt{2} +o_P(1).
\]
\end{theo}

Thus, under the assumptions from part~(ii), one has
\[
\overline{B^{H_k,\tau_k}} -  \sqrt{\ln n_k} \todi \zeta_0/\sqrt{2} \quad\mbox{as}\quad k\to\infty.
\]
Note also that the conditions $H_k\ln (n_k\delta_k) \to 0$ and $H_k \ln\delta_k\to0$ from parts~(i) and~(ii), respectively, are automatically met in the case of ``uniform grids" $\tau_k$ (when $\delta (\tau_k)=1/n_k$).

\begin{figure}[ht]
\centering
\includegraphics[width=7cm]{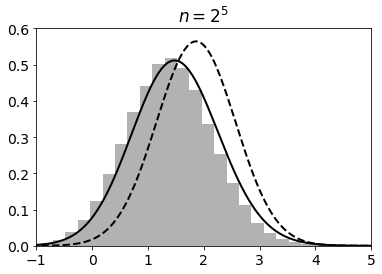}\hfill
\includegraphics[width=7cm]{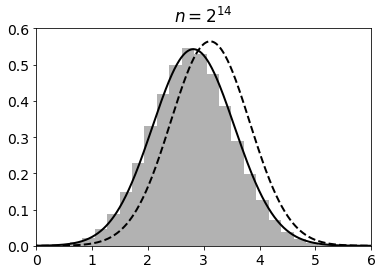}

\caption{The histograms show the empirical distributions of
$\overline{B^{H,\tau}}$ for $10^5$ simulated paths of the fBM $B^H$ with the uniform partition
$\tau = \{i/n\}_{1\le i \le n}$ and $H = (\ln n)^{-2}$. The dashed lines show
the approximating normal densities, and the solid lines the approximations
by the convolutions of the scaled Gumbel and normal densities.}
\label{fig1}
\end{figure}

Simulations indicate that in fact, in accordance with~\eqref{toGu}, a better
approximation to the law of $\overline{B^{H_k,\tau_k}}$ is given by the distribution $D_n(x): = (\Lambda_n * \Phi)(\sqrt2 x),$ the convolution being that  of the scaled version of the Gumbel law
$\Lambda_n(x) = \Lambda(a_n(x-b_n))$ with  the standard  normal distribution. The curves in Fig.~\ref{fig1} are the fitting normal density (dashed lines) and the density of $D_n$ (solid  lines), where $a_n, b_n$ were chosen according to~\eqref{ab}, overlayed upon the histograms constructed from the respective simulations. However, establishing the validity of that second order approximation analytically is much
harder than the analysis in the present note and may require more refined
techniques.

\section{The proof of the theorem}

(i)~Let $W$ be a standard BM process independent
of $\{\zeta_t\}$. Set $s_{k,i}:= (t_{k,i})^{2H_k}$, $i=1,\ldots,n_k$, and
introduce   random vectors $X^k, Y^k\in\R^{n_k}$ with the respective components
\[
X^k_i:= (s_{k,i}^{1/2}\zeta_i - W_{s_{k,i}})/\sqrt{2}, \qquad Y^k_i:= (\zeta_i - W_{s_{k,i}})/\sqrt2.
\]
First we  show that
\begin{equation}
\overline{B^{H_k,\tau_k}}\lest \overline{X^k},
\label{XYineq}
\end{equation}
then give an upper bound for $\overline{X^k}$ in terms of $\overline{Y^k},$ and finally demonstrate that that bound is of the form of the RHS of~\eqref{<=}.

Clearly, $\exn X^k=0$ and
\begin{align*}
\Co (X^k_i, X^k_j)
& = 2^{-1}\bigl(s_{k,i}^{1/2}s_{k,j}^{1/2}\Co (\zeta_i, \zeta_j) + \Co (W_{s_{k,i}}, W_{s_{k,j}})\bigr)
\\
&
=2^{-1}\bigl(s_{k,i} \delta_{ij} + s_{k,i}\wedge s_{k,j})\bigr), \quad 1\le i,j\le n_k,
\end{align*}
where $\delta_{ij}$ is Kronecker's delta. Therefore,
\begin{align}
\exn X^k_i = \exn B^{H_k,\tau_k}_i, \quad \Va  X^k_i = \Va B^{H_k,\tau_k}_i, \quad 1\le i\le n_k,
\label{VarVar}
\end{align}
and, for $1\le i< j \le n_k,$ one has
\begin{equation}
\begin{split}
\Co (X^k_i, X^k_j)&=\mbox{$\frac12$} s_{k,i}
  < \mbox{$\frac12$}\bigl( s_{k.i} + s_{k,j} -s_{k,j}
(1- t_{k,i}/t_{k,j}  )^{2H_k}\bigr)
 \\&
 = \Co (B^{H_k,\tau_k}_i, B^{H_k,\tau_k}_j).
\end{split}
\label{CoCo}
\end{equation}
Now  \eqref{XYineq} immediately follows from Slepian's lemma~\cite{Sl}.

Next let $i(k): = \mathrm{argmax}_{1\le i \le n_k} X^k_i$, which is clearly well-defined a.s. Since
$s_{k,i}\le 1,$ it is easy to see that 
\begin{equation}
\overline{X^k} 
 \le \overline{Y^k} \ind(\zeta_{i(k)} \ge 0)  - 2^{-1/2}
W_{s_{k,i(k) }} \ind(\zeta_{i(k)} < 0).
\label{X^k bound}
\end{equation}
We will now show that 
\begin{equation}
\overline{Y^k}\le \sqrt{\ln n_k }-W_1/\sqrt{2}+o_P(1).
\label{Y^k bound}
\end{equation}
The assumption that $H_k(\ln n_k)^{1/2}\to 0$ ensures that it is possible to
choose a sequence $\ep_k>0$ such that the following
relations hold as $k\to\infty$: 
\begin{align}
&   \ep_k\to 0, \qquad m_k:= \ep_k n_k\in \N, 
\notag
\\ 
&\frac{|\ln \ep_k|}{\ln n_k} \to0,\qquad
\frac{|\ln \ep_k|}{\sqrt{\ln n_k}}\to \infty,
\label{ep_1}
 \\
&m_k  \to\infty,\qquad
H_k|\ln \ep_k|  \to0.
\label{ep_2}
\end{align}
Indeed, one can set $\ep_k:= e^{-N_k\sqrt{\ln n_k}}$ with a quantity
$N_k\to\infty$ such that
$ N_k (\ln n_k)^{1/2} = o \bigl(H^{-1}_k\wedge \ln n_k\bigr) $ (for example,
$N_k := (H_k(\ln n_k)^{1/2})^{-1/2} \wedge (\ln n_k)^{1/4}$, adjusted if necessary to ensure that $m_k\in\N$).

Now set $C_{k,1}:=\{i: 1\le i \le m_k \}$,
$C_{k,2}:=\{i:  m_k < i \le n_k\}$  and let
\[
M_{k,j}:= \max_{i\in C_{k,j}} \bigl( \zeta_i - W_{s_{k,i}}\bigr) , \quad
j=1,2,
\]
so that $ \overline{Y^k}= (M_{k,1} \vee M_{k,2})/\sqrt2$.

To bound $M_{k,1},$  note that
\[
x_k:=\sqrt{2\ln m_k}= \sqrt{2\ln n_k \Bigl(1 + \frac{\ln \ep_k}{\ln n_k}\Bigr)}
 \le \sqrt{2\ln n_k}\Bigl(1 + \frac{\ln \ep_k}{2\ln n_k}\Bigr)
 = \sqrt{2\ln n_k} - 2h_k,
\]
where in view of~\eqref{ep_1} one has 
\begin{equation}
h_k:= |\ln \ep_k |/(2\sqrt{2\ln n_k})\to \infty.
 \label{h_n}
\end{equation}
Using the standard  Mills' ratio bound for the normal distribution, we have
\begin{equation}
\pr (\overline{\zeta^{m_k}} >x_k ) \le m_k \pr (\zeta_1 >x_k) \le \frac{m_k e^{-x^2_k/2}}{\sqrt{2\pi }x_k}
= \frac1{\sqrt{4\pi \ln m_k}}
\to0
\label{zeta_m}
\end{equation}
in view of~\eqref{ep_2}. Setting $\underline{W_1}:=\min_{0\le t\le 1} W_t$, we obtain that
\begin{align}
\pr (M_{k,1} > \sqrt{2\ln n_k} -  h_k)
   & \le \pr (\overline{\zeta^{m_k}} - \underline{W_1} > \sqrt{2\ln n_k} -  h_k)
  \notag
  \\
  & \le \pr (\overline{\zeta^{m_k}}   > \sqrt{2\ln n_k} -  2h_k)+ \pr (  - \underline{W_1} > h_k)\to0
\label{M1}
\end{align}
by~\eqref{h_n} and~\eqref{zeta_m}.

Now we turn to the term~$M_{k,2}$. As $W$ has continuous trajectories, there
exist $\theta_k\in [s_{k,m_k},1]$, which depend  on the trajectory of~$W$,  such that
\begin{align}
M_{k,2} = \max_{m_k < i\le n_k}\zeta_i - W_{\theta_k} \le \overline{\zeta^{n_k}} -W_1 + o_P(1),
\label{M2}
\end{align}
where the last relation  holds as  $W_{\theta_k}\to W_1$ because $\theta_k\to1$
since 
\begin{align}
s_{k,m_k} \ge (m_k\delta_k )^{2H_k} = \ep_k^{2H_k} (n_k\delta_k )^{2H_k}\to1 
\label{s_bound}
\end{align}
due to  the assumption that $H_k \ln (n_k\delta_k )\to 0$  and~\eqref{ep_2}.

Since $\overline{\zeta^{n_k}}=\sqrt{2\ln n_k }+o_P(1) $ in view
of~\eqref{toGu}, from \eqref{M1} and \eqref{M2} we obtain that
$M_{k,1}\vee M_{k,2}\le \sqrt{2\ln n_k} - W_1 + o_P(1)$, which proves
\eqref{Y^k bound}.

Now  observe
that obviously 
\[
 -   W_{s_{k,i(k) }}\le \sqrt{2\ln n_k} - W_1 + o_P(1) 
\]
and $W_1\deq -\zeta_0$. That, together with   \eqref{XYineq}, \eqref{X^k bound} and~\eqref{Y^k bound}, completes the proof of  part~(i) of the theorem.

\medskip
(ii)~Consider the differences
\[
d_{k,ij}:= \bigl(\Co B^{H_k,\tau_k}- \Co X^k\bigr)_{ij}\ge 0, \quad 1\le i,j \le n_k
\]
(cf.~\eqref{VarVar}, \eqref{CoCo}). Note that $d_{k,ii}=0,$ $1\le i\le n_k,$ by~\eqref{VarVar}, and that for $i<j$ one has
\begin{align*}
d_{k,ij}  & = \frac12 \biggl[\Bigl(\frac{j}{n_k}\Bigr)^{2H_k}
 - \Bigl(\frac{j-i}{n_k}\Bigr)^{2H_k}\biggr]
  \le
  \frac12 \biggl[1
 - \Bigl(\frac{1}{n_k}\Bigr)^{2H_k}\biggr] \le  H_k\ln n_k := q_k
\end{align*}
since $1-1/x \le \ln x$ for all $x>0$. Denoting by $I_k:=(\delta_{ij}) $ and $J_k:=(1)$ the unit and all-ones $(n_k\times n_k)$-matrices, respectively, we conclude that
\begin{align}
(\Co B^{H_k,\tau_k}+ q_k I_k )_{ij}\le ( \Co X^k + q_k J_k )_{ij}, \quad 1\le i, j\le n_k,
\label{CoCo1}
\end{align}
with equalities holding for $i=j$.

On the LHS of~\eqref{CoCo1} we have got the entries of the covariance matrix of the random vector $ B^{H_k,\tau_k}+ q_k^{1/2}\zeta^{n_k} $ (assuming that $\{\zeta_t\}$ is independent of~$ B^{H_k}$), whereas on the RHS are those for the vector $X^k+ q_k^{1/2}\zeta_0$ (addition with a scalar is understood in the component-wise sense). Since the means of those random vectors are zeros, by Slepian's lemma one has
\[
\overline{ B^{H_k,\tau_k}+ q_k^{1/2}\zeta^{n_k}} \gest \overline{X^k+ q_k^{1/2}\zeta_0}
 = \overline{X^k}+ q_k^{1/2}\zeta_0  = \overline {X^k} + o_P(1).
\]
Using~\eqref{toGu}, we have
\[
q_k^{1/2} \overline{\zeta^{n_k}}=q_k^{1/2}\sqrt{2\ln n_k} +o_P(1) = o_P(1)
\]
as $q_k \ln n_k= H_k(\ln n_k)^2=o(1)$ by assumption. Hence, by the lemma from the Appendix, one has 
\begin{equation}
\overline{ B^{H_k,\tau_k} } \ge \overline{ B^{H_k,\tau_k}+ q_k^{1/2}\zeta^{n_k}}
- q_k^{1/2}\overline{\zeta^{n_k}}
\gest  \overline{X^k}+ o_P(1).
\label{aaa}
\end{equation}

On the event $A_k = \{\max_{m_k< i \le n_k}\zeta_i \ge 0\} $ we have
\[
2^{1/2} \overline{X^k}
 \ge \max_{m_k < i\le n_k}\bigl( s_{k,i}^{1/2}\zeta_i - W_{s_{k,i}}\bigr)
  \ge s_{k,m_k}^{1/2} \max_{m_k <i\le n_k} \zeta_i 
  + 
  \min_{s_{k,m_k} \le t \le 1 }W_t.
\]
In view of the first two relations in~\eqref{s_bound}, the second relation in~\eqref{ep_2} and the assumption of part~(ii) of the theorem, we have
$s_{k,m_k} \to 1$ as $k\to\infty.$ Therefore, 
\[
\min_{s_{k,m_k}  \le t \le 1 }W_t \deq \zeta_0 +o_P(1).
\] 
Since clearly $\pr (A_k)\to 1,$ we obtain that
\[
2^{1/2} \overline{X^k}
  \gest s_{k,m_k}^{1/2} \max_{m_k < i\le n_k} \zeta_i +\zeta_0 + o_P(1).
\]
For the first term on the RHS, using~\eqref{toGu}, one has
\[
\max_{m_k < i\le n_k} \zeta_i \deq\overline{\zeta^{(1-\ep_k )n_k}}
=
\sqrt{2\ln ((1-\ep_k)n_k)} +o_P(1)= \sqrt{2\ln n_k } +o_P(1)
\]
as clearly $\ep_k\sqrt{\ln n_k} =o(1)$. Thus,
$\overline{X^k} \gest \sqrt{\ln n_k} + \zeta_0/\sqrt{2} + o_P(1)$.  To complete the proof of part~(ii) of the theorem, it remains to combine the last bound with~\eqref{aaa} and  again use the lemma from the Appendix.

\appendix
\section*{Appendix}
The following simple lemma was used in the proof of the theorem.

\begin{lemo}
Suppose $X,Y$ are two random variables such that $X$ has a continuous
distribution and $X\gest Y$, while $Z$ is a random variable defined on the
same probability space as $X$. Then there exist random variables $Y', Z'$
such that $X+Z \gest Y' + Z'$ and $Y\deq Y'$, $Z\deq Z'$.
\end{lemo}

In particular, if $X_n \gest Y_n$ and $Z_n \topr0$ as $n\to\infty$, then
$X_n+Z_n \gest Y'_n + o_P(1)$, where $Y'_n\deq Y_n$ for all $n$. In fact, the assumption
that $X$ has a continuous distribution can be relaxed, by that is not
necessary for us.   

Note that if $X,Y,Z$ are defined on the same probability space, then the
inequality $X\gest Y$ does not necessarily imply that $X+Z \gest Y+Z$. Here
is a counterexample: let $X$ be a uniform random variable on $[0,1]$ and set $Y:=Z:=1-X$.

The proof of the lemma readily follows from the explicit construction
$Y': = F^{(-1)}_Y(F_X(X))$, $Z':=Z$, where $F_X$, $F_Y$ denote the corresponding
distribution functions, $F^{(-1)}_Y$ the generalized inverse of~$F_Y$. Then $X,Y',Z$ are defined on the same probability
space, and $X+Z \ge Y'+Z'$ with probability one.

\section*{Acknowledgments}
K.~Borovkov was supported by the ARC  Discovery grant DP150102758. Some of the presented results   were  obtained when
the authors were taking part in the   ``Mathematics of Risk'' program ran at the MATRIX Research Institute in~2017. We are grateful to MATRIX for
its support and hospitality.

\setlength{\bibsep}{0.5em}

\end{document}